\documentclass[a4paper,9pt]{article}
\usepackage{amsmath, amsfonts, amssymb}
\usepackage[english, russian]{babel}
\usepackage[cp1251]{inputenc}

\usepackage{graphicx}
\begin{document}
\sloppy
\begin{center}
\textbf{Fundamental solutions  of the generalized Helmholtz equation with several singular coefficients and confluent hypergeometric functions of many variables}\\[5pt]
\textbf{Ergashev T.G.\\}
\medskip
{ Institute of Mathematics, Uzbek
Academy of Sciences,  Tashkent, Uzbekistan. \\

{\verb ergashev.tukhtasin@gmail.com }}\\
\end{center}

\begin{quote} 

In this paper, we introduce a new class of confluent hypergeometric functions of many variables, study their properties, and determine a system of partial differential equations that this function satisfies. It turns out that all the fundamental solutions of the generalized Helmholtz equation with several singular coefficients are written out through the newly introduced confluent hypergeometric function. Using the expansion formula established here for the confluent function, the order of the singularity of the fundamental solutions of the elliptic equation under this consideration is determined.

\end{quote}

\section{Introduction and Preliminaries}

An investigation of applied problems related to heat conduction and dynamics, electromagnetic oscillations and aerodynamics, quantum mechanics and potential theory leads to the study of various hypergeometric functions. The great success of the theory of hypergeometric functions in one variable has stimulated the development of corresponding theory in two and more variables. Books \cite{E1953} and \cite{SK} are respectively devoted to a systematic presentation of the results on the hypergeometric functions of two and three variables. In the theory of hypergeometric functions, an increase in a number of variables will always be accompanied by a complication in the study of the function of several variables. Therefore, the decomposition formulas that allow us to represent the hypergeometric function of several variables through an infinite sum of products of several hypergeometric functions in one variable are very important, and this, in turn, facilitate the process of studying the properties of multidimensional functions. In the literature, hypergeometric functions are divided  into two types: complete and confluent (for definitions see \cite{E1953}). In all respects, confluent hypergeometric functions including the decomposition formulas, have been little studied in comparison with other types of hypergeometric functions, especially when the dimension of the variables exceeds two. We only note the works \cite{{J},{E}}, in which some confluent hypergeometric functions of three variables were сonsidered. In this paper we define a new class of confluent hypergeometric functions of several variables, study their properties, and establish decomposition formulas. An important application of confluent functions was found. It turns out that all fundamental solutions of the generalized Helmholtz equation with singular coefficients are written out through one new introduced confluent hypergeometric function of several variables. Using the decomposition formulas, the order of the singularity of the fundamental solutions of the elliptic equation which mentioned above is determined.

    Let us consider the generalized Helmholtz equation with a several singular coefficients
\numberwithin{equation}{section}
\begin{equation} \label{eq11}
\sum\limits_{i=1}^{m}{\frac{{{\partial }^{2}}u}{\partial x_{i}^{2}}}+\sum\limits_{j=1}^{n}{\frac{2{{\alpha }_{j}}}{{{x}_{j}}}\frac{\partial u}{\partial {{x}_{j}}}}-{{\lambda }^{2}}u=0
\end{equation}
in the domain $R_{m}^{n+}:=\left\{ \left( {{x}_{1}},{{x}_{2}},...,{{x}_{m}} \right):{{x}_{1}}>0,{{x}_{2}}>0,...,{{x}_{n}}>0 \right\},$ where $m\ge 2$ is a dimension of the Euclidean space; $n\ge 1$ is a number of the singular coefficients;  $m\ge n;$ ${{\alpha }_{j}}$  are real constants and $0<2{{\alpha }_{j}}<1$$(j=\overline{1,n});$ $\lambda $ is real or pure imaginary constant.

Various modifications of equation (\ref{eq11}) in the two- and three-dimensional cases were considered in many papers \cite{{H07},{HK},{UK},{U}}. However, relatively few papers have been devoted to finding the fundamental solutions when the dimension of equation exceeds three, we only note the works \cite{{U},{Garip},{EH},{E19},{UE}}. For example, in the paper \cite{UE}, fundamental solutions of equation (\ref{eq11}) when a dimension of the Euclidean space equals to the number of the singular coefficients ($m=n$) are found and investigated. As it is known, fundamental solutions play an important role in the study of partial differential equations. The formulation and solution of many local and nonlocal boundary value problems are based on these solutions. It turns out that all the fundamental solutions of equation (\ref{eq11}) are written out through the confluent hypergeometric function of the $\left( n+1 \right)$ variables. In this paper, we find all the fundamental solutions of equation (\ref{eq11}), write out the expansion formulas for them. By using these formulas we establish that the fundamental solutions of equation (\ref{eq11}) have a singularity of order $1/{{r}^{m-2}}$ at $r\to 0$.

We give some following notation and definitions, which will be used in the next sections.

A symbol ${{\left( \text{ }\!\!\kappa\!\!\text{ } \right)}_{\nu }}$ denotes the general Pochhammer symbol or the shifted factorial, since ${{\left( 1 \right)}_{l}}=l!$ $\left( l\in \mathbb{N}\bigcup \left\{ 0 \right\};\mathbb{N}:=\left\{ 1,2,... \right\} \right),$ which is defined (for $\text{ }\!\!\kappa\!\!\text{ }\text{,}\,\nu \in \mathbb{C}$), in terms of the familiar Gamma function, by
$${{\left( \text{ }\!\!\kappa\!\!\text{ } \right)}_{\text{ }\!\!\nu\!\!\text{ }}}:=\frac{\Gamma \left( \text{ }\!\!\kappa\!\!\text{ }+\text{ }\!\!\nu\!\!\text{ } \right)}{\Gamma \left( \text{ }\!\!\kappa\!\!\text{ } \right)}=\left\{ \begin{array}{*{20}c}
   1 & \left( \nu =0;\,\kappa \in \mathbb{C}\backslash \left\{ 0 \right\} \right),  \\
   \kappa \left( \kappa +1 \right)...\left( \kappa +l-1 \right) & \left( \nu =l\in \mathbb{N};\,\kappa \in \mathbb{C} \right),  \\
\end{array} \right.$$
it is being understood conventionally that  ${{\left( 0 \right)}_{0}}:=1$ assumed tacitly that the $\Gamma -$quotient exists.

We introduce the following functions:
$$F\left( a,b;c;x \right)=\sum\limits_{m=0}^{\infty }\frac{{{\left( a \right)}_{m}}{{\left( b \right)}_{m}}}{{{\left( c \right)}_{m}}m!}{{x}^{m}},\,|x|<1;$$
$${{H}_{2}}\left( a,{{b}_{1}},{{d}_{1}},{{d}_{2}};{{c}_{1}};x,y \right)=\sum\limits_{m,n=0}^{\infty }{\frac{{{\left( a \right)}_{m-n}}{{\left( {{b}_{1}} \right)}_{m}}{{\left( {{d}_{1}} \right)}_{n}}{{\left( {{d}_{2}} \right)}_{n}}}{m!n!{{\left( {{c}_{1}} \right)}_{m}}}{{x}^{m}}{{y}^{n}}},|x|<1,\,\,|y|<1/\left(1+|x| \right);$$
\begin{equation} \label{eq12}
F_{A}^{(n)}\left( a,{{b}_{1}},...,{{b}_{n}};{{c}_{1}},...,{{c}_{n}};{{x}_{1}},...,{{x}_{n}} \right)=\sum\limits_{{{m}_{1}},...{{m}_{n}}=0}^{\infty }{\frac{{{\left( a \right)}_{{{m}_{1}}+...+{{m}_{n}}}}{{\left( {{b}_{1}} \right)}_{{{m}_{1}}}}...{{\left( {{b}_{n}} \right)}_{{{m}_{n}}}}}{{{m}_{1}}!...{{m}_{n}}!{{\left( {{c}_{1}} \right)}_{{{m}_{1}}}}...{{\left( {{c}_{n}} \right)}_{{{m}_{n}}}}}x_{1}^{{{m}_{1}}}...x_{n}^{{{m}_{n}}}},
\end{equation}
 $$\,\left| {{x}_{1}} \right|+\left| {{x}_{2}} \right|+...+\left| {{x}_{n}} \right|<1;$$
 $$F_{B}^{(n)}\left( {{a}_{1}},...,{{a}_{n}},{{b}_{1}},...,{{b}_{n}};c;{{x}_{1}},...,{{x}_{n}} \right)=\sum\limits_{{{m}_{1}},...{{m}_{n}}=0}^{\infty }{\frac{{{\left( {{a}_{1}} \right)}_{{{m}_{1}}}}...{{\left( {{a}_{n}} \right)}_{{{m}_{n}}}}{{\left( {{b}_{1}} \right)}_{{{m}_{1}}}}...{{\left( {{b}_{n}} \right)}_{{{m}_{n}}}}}{{{m}_{1}}!...{{m}_{n}}!{{\left( c \right)}_{{{m}_{1}}+...+{{m}_{n}}}}}x_{1}^{{{m}_{1}}}...x_{n}^{{{m}_{n}}}},$$ $$\max \left( \,\left| {{x}_{1}} \right|,\left| {{x}_{2}} \right|,...,\left| {{x}_{n}} \right| \right)<1.$$

Here $F$ is a hypergeometric function of Gauss \cite[p.56, (2)]{E1953}; $H_2$ is Horn's function in two variables \cite[p.225, (14)]{E1953};  $F_{A}^{(n)}$ and $F_{B}^{(n)}$ are Lauricella's hypergeometric functions $n\in \mathbb{N}$ variables \cite{L} (see also \cite[p.33, (1)-(2)]{SK}).

\section{Generalizations of Lauricella functions and its confluent forms}

An interesting unification (and generalization) of Lauricella’s multivariable functions  and Horn’s two variables functions     was considered by Erdelyi in 1939, who defined his general functions in the form \cite{E1939} (see also \cite[p.74, (4b)]{SK}):
     $${{H}_{n+p,n}}\left( a,{{b}_{1}},...,{{b}_{n+p}},{{d}_{n+1}},...,{{d}_{n+p}};{{c}_{1}},...,{{c}_{n}};{{\xi }_{1}},...,{{\xi }_{n}},{{\eta }_{1}},...,{{\eta }_{p}} \right)=$$
$$=\sum\limits_{{{m}_{1}},...,{{m}_{n+p}}=0}^{\infty }{{}}\frac{{{\left( a \right)}_{{{m}_{1}}+...+{{m}_{n}}-{{m}_{n+1}}-...-{{m}_{n+p}}}}{{\left( {{b}_{1}} \right)}_{{{m}_{1}}}}...{{\left( {{b}_{n+p}} \right)}_{{{m}_{n+p}}}}{{\left( {{d}_{n+1}} \right)}_{{{m}_{n+1}}}}...{{\left( {{d}_{n+p}} \right)}_{{{m}_{n+p}}}}}{{{m}_{1}}!...{{m}_{n+p}}!{{\left( {{c}_{1}} \right)}_{{{m}_{1}}}}...{{\left( {{c}_{n}} \right)}_{{{m}_{n}}}}}$$
\begin{equation} \label{eq21}
\cdot\text{ }\!\!\xi\!\!\text{ }_{1}^{{{m}_{1}}}...\text{ }\!\!\xi\!\!\text{ }_{n}^{{{m}_{n}}}\eta _{1}^{{{m}_{n+1}}}...\eta _{p}^{{{m}_{n+p}}},
\end{equation}
where  $p$ and  $n$ are nonnegative integers.

Evidently, we have
$${{H}_{n,n}}=F_{A}^{(n)},{{H}_{n,0}}=F_{B}^{(n)},{{H}_{2,1}}={{H}_{2}}.$$

From the hypergeometric function (\ref{eq21}) we shall define the following confluent hypergeometric function
$${\rm{H}}_{A}^{(n,p)}\left( a,{{b}_{1}},...,{{b}_{n}};{{c}_{1}},...,{{c}_{n}};{{\xi }_{1}},...,{{\xi }_{n}},{{\eta }_{1}},...,{{\eta }_{p}} \right)$$
$$=\underset{\varepsilon \to 0}{\mathop{\lim }}\,{{H}_{n+p,n}}\left( a,{{b}_{1}},...,{{b}_{n}},\underbrace{\frac{1}{\varepsilon },...,\frac{1}{\varepsilon }}_{2p};{{c}_{1}},...,{{c}_{n}};{{\xi }_{1}},...,{{\xi }_{n}},{{\text{ }\!\!\varepsilon\!\!\text{ }}^{2}}{{\eta }_{1}},...,{{\text{ }\!\!\varepsilon\!\!\text{ }}^{2}}{{\eta }_{p}} \right).$$

For the determination of the hypergeometric function ${\rm{H}}_{A}^{(n,p)}$  the equality \cite[p.124]{App} $\underset{\text{ }\!\!\varepsilon\!\!\text{ }\to 0}{\mathop{\lim }}\,{{\left( 1/\text{ }\!\!\varepsilon\!\!\text{ } \right)}_{q}}\cdot {{\text{ }\!\!\varepsilon\!\!\text{ }}^{q}}=1$ ($q$ is a natural number) has been used and the found confluent hypergeometric function is represented as
\begin{equation*}
{\rm{H}}_{A}^{(n,p)}\left( a,{{b}_{1}},...,{{b}_{n}};{{c}_{1}},...,{{c}_{n}};\xi ,\eta  \right)={\rm{H}}_{A}^{(n,p)}\left[ \begin{matrix}
   a,{{b}_{1}},...,{{b}_{n}};  \\
   {{c}_{1}},...,{{c}_{n}};  \\
\end{matrix}\xi \text{,}\eta  \right]
\end{equation*}
\begin{equation} \label{eq22}
=\sum\limits_{{{m}_{1}},...,{{m}_{n+p}}=0}^{\infty }{\frac{{{\left( a \right)}_{{{m}_{1}}+...+{{m}_{n}}-{{m}_{n+1}}-...-{{m}_{n+p}}}}{{\left( {{b}_{1}} \right)}_{{{m}_{1}}}}...{{\left( {{b}_{n}} \right)}_{{{m}_{n}}}}}{{{m}_{1}}!...{{m}_{n+p}}!{{\left( {{c}_{1}} \right)}_{{{m}_{1}}}}...{{\left( {{c}_{n}} \right)}_{{{m}_{n}}}}}\xi _{1}^{{{m}_{1}}}...\xi _{n}^{{{m}_{n}}}\eta _{1}^{{{m}_{n+1}}}...\eta _{p}^{{{m}_{n+p}}}},
\end{equation}
where $\xi :=\left( {{\xi }_{1}},...,{{\xi }_{n}} \right),\,\,\eta :=\left( {{\eta }_{1}},...,{{\eta }_{p}} \right)$ and $\,\left| {{\xi }_{1}} \right|+...+\,\left| {{\xi }_{n}} \right|<1.$

The particular cases of confluent hypergeometric functions were known: in the case of two variables \cite[p.231, (31)]{E1953}
$${{\rm{H} }_{3}}\left( a,b;c;x,y \right)=\sum\limits_{m,n=0}^{\infty }{{}}\frac{{{\left( a \right)}_{m-n}}{{\left( b \right)}_{m}}}{m!n!{{\left( c \right)}_{m}}}{{x}^{m}}{{y}^{n}},\,\,\left| x \right|<1$$
and in the case of three variables \cite{H07}
$${{\rm{A} }_{2}}\left( a,{{b}_{1}},{{b}_{2}};{{c}_{1}},{{c}_{2}};x,y,z \right)=\sum\limits_{m,n,k=0}^{\infty }{{}}\frac{{{\left( a \right)}_{m+n-k}}{{\left( {{b}_{1}} \right)}_{m}}{{\left( {{b}_{2}} \right)}_{n}}}{m!n!k!{{\left( {{c}_{1}} \right)}_{m}}{{\left( {{c}_{2}} \right)}_{n}}}{{x}^{m}}{{y}^{n}}{{z}^{k}},\,\,\left| x \right|+\left| y \right|<1.$$

The confluent hypergeometric function ${\rm{H}}_{A}^{(n,p)}$ has the following formula of derivation:
\begin{equation*}
\frac{{{\partial }^{{{i}_{1}}+...+{{i}_{n}}+{{j}_{1}}+...+{{j}_{p}}}}}{\partial \xi _{1}^{{{i}_{1}}}...\partial \xi _{n}^{{{i}_{n}}}\partial \eta _{1}^{{{j}_{1}}}...\partial \eta _{p}^{{{j}_{p}}}}{\rm{H}}_{A}^{(n,p)}\left( a,{{b}_{1}},...,{{b}_{n}};{{c}_{1}},...,{{c}_{n}};\xi ,\eta  \right)=\frac{{{\left( a \right)}_{{{i}_{1}}+...+{{i}_{n}}-{{j}_{1}}-...-{{j}_{p}}}}{{\left( {{b}_{1}} \right)}_{{{i}_{1}}}}...{{\left( {{b}_{n}} \right)}_{{{i}_{n}}}}}{{{\left( {{c}_{1}} \right)}_{{{i}_{1}}}}...{{\left( {{c}_{n}} \right)}_{{{i}_{n}}}}}
\end{equation*}
\begin{equation} \label{eq23}
\cdot {\rm{H}}_{A}^{(n,p)}\left( a+{{i}_{1}}+...+{{i}_{n}}-{{j}_{1}}-...-{{j}_{p}},{{b}_{1}}+{{i}_{1}},...,{{b}_{n}}+{{i}_{n}};{{c}_{1}}+{{i}_{1}},...,{{c}_{n}}+{{i}_{n}};\xi ,\eta  \right).
\end{equation}

Using the formula of derivation (\ref{eq23}) it is easy to show that the confluent hypergeometric function in (\ref{eq22}) satisfies  the following system of hypergeometric equations
\begin{equation} \label{eq24}
\left\{ \begin{array}{*{20}c}
   {{\xi }_{i}}\left( 1-{{\xi }_{i}} \right){{\omega }_{{{\xi }_{i}}{{\xi }_{i}}}}-{{\xi }_{i}}\sum\limits_{j=1,j\ne i}^{n}{{{\xi}_{j}}{{\omega }_{{{\xi }_{i}}{{\xi }_{j}}}}}+{{\xi }_{i}}\sum\limits_{j=1}^{p}{{{\eta }_{j}}{{\omega }_{{{\xi }_{i}}{{\eta}_{j}}}}}+\left[ {{c}_{i}}-\left( a+{{b}_{i}}+1 \right){{\xi }_{i}} \right]{{\omega }_{{{\xi }_{i}}}}\\
   -{{b}_{i}}\sum\limits_{j=1,j\ne i}^{n}{{{\xi }_{j}}{{\omega }_{{{\xi }_{j}}}}}+{{b}_{i}}\sum\limits_{j=1}^{p}{{{\eta }_{j}}{{\omega}_{{{\eta }_{j}}}}}-a{{b}_{i}}\omega =0,\,\,\,i=\overline{1,n},  \\
   \sum\limits_{l=1}^{p}{{{\eta }_{l}}{{\omega }_{{{\eta }_{l}}{{\eta }_{j}}}}}-\sum\limits_{l=1}^{n}{{{\xi }_{l}}{{\omega }_{{{\xi}_{l}}{{\eta }_{j}}}}}+\left( 1-a \right){{\omega }_{{{\eta }_{j}}}}+\omega =0,\,\,j=\overline{1,p},  \\
\end{array} \right.
\end{equation}
where $\omega \left( \xi ,\eta  \right)={\rm{H}} _{A}^{(n,p)}\left( a,{{b}_{1}},...,{{b}_{n}};{{c}_{1}},...,{{c}_{n}};\xi ,\eta  \right).$

We note that particular cases (i.e. $n=1,2,3$ for $p=1$) of the system (\ref{eq24}) are found in \cite{{E1953},{U},{Garip},{EH},{E19}}.

Having substituted  $\omega \left( \xi ,\eta  \right)=\xi _{1}^{{{\tau }_{1}}}...\xi _{n}^{{{\tau }_{n}}}\eta _{1}^{{{\nu }_{1}}}...\eta _{p}^{{{\nu }_{p}}}\psi \left( \xi ,\eta  \right)$ in the system (\ref{eq24}), it is possible to find ${{2}^{n}}$ linearly independent solutions of system (\ref{eq24}), which are given as follows (for details, see \cite{UE}):
$$\left. {\rm{H}}_{A}^{(n,p)}\left( a,{{b}_{1}},...,{{b}_{n}};{{c}_{1}},...,{{c}_{n}};\xi ,\eta  \right) \right\}C_{n}^{0}=1$$
$$\left. \begin{array}{*{20}c}
   \xi _{1}^{1-{{c}_{1}}}{\rm{H}}_{A}^{(n,p)}\left( a,1+{{b}_{1}}-{{c}_{1}},{{b}_{2}},...,{{b}_{n}};2-{{c}_{1}},{{c}_{2}},...,{{c}_{n}};\xi ,\eta  \right)  \\
   ...................................................................  \\
   \xi _{n}^{1-{{c}_{n}}}{\rm{H}} _{A}^{(n,p)}\left( a,{{b}_{1}},...,{{b}_{n-1}},1+{{b}_{n}}-{{c}_{n}};{{c}_{1}},...,{{c}_{n-1}},2-{{c}_{n}};\xi ,\eta  \right)  \\
\end{array} \right\}C_{n}^{1}=n$$
$$.............................................................$$
$$\left. \xi _{1}^{1-{{c}_{1}}}...\xi _{n}^{1-{{c}_{n}}}{\rm{H}}_{A}^{(n,p)}\left( a,1+{{b}_{1}}-{{c}_{1}},...,1+{{b}_{n}}-{{c}_{n}};2-{{c}_{1}},...,2-{{c}_{n}};\xi ,\eta  \right) \right\}C_{n}^{n}=1.$$

Here,  $C_{n}^{i}=\displaystyle\frac{n!}{i!\left( n-i \right)!}$  and  it is easy to see that $1+C_{n}^{1}+C_{n}^{2}+...+C_{n}^{n-1}+1={{2}^{n}}.$

By virtue of symmetry of the function ${\rm{H}}_{A}^{(n,p)}$ with respect to the parameters ${{b}_{1}},$…,${{b}_{n}},$${{c}_{1}},$…,${{c}_{n}}$, it is possible to group the above linearly independent solutions of the system of hypergeometric equations (\ref{eq24}). As a result, a number of solutions of the system (\ref{eq24}), which are necessary to further studies, will decrease. Thus, all solutions of the system (\ref{eq24}) are expressed by the formula:
\begin{equation} \label{eq25}
{{\omega }_{k}}\left( \xi ,\eta  \right)={{C}_{k}}\prod\limits_{i=1}^{k}{\left[ \xi _{i}^{1-{{c}_{i}}} \right]}\cdot {\rm{H}}_{A}^{\left( n,p \right)}\left[ \begin{array}{*{20}c} a,{{b}_{1}}+1-{{c}_{1}},...,{{b}_{k}}+1-{{c}_{k}},{{b}_{k+1}},...,{{b}_{n}};\\
2-{{c}_{1}},...,2-{{c}_{k}},{{c}_{k+1}},...,{{c}_{n}};\\ \end{array}
\xi,\,\eta \right], \,k=\overline{0,n},
\end{equation}
where  ${{C}_{k}}$ are arbitrary constants.

\section{Decomposition formulas}

For a given multivariable function, it is useful to find a decomposition formula which would express the multivariable function in terms of products of several simpler hypergeometric functions involving fewer variables. For this purpose, Burchnall and Chaundy \cite{{BC1},{BC2}} found a number of expansions of double hypergeometric functions in series of simpler hypergeometric functions. Their method is based on the inverse pair of the symbolic operators
\begin{equation} \label{eq31}
\nabla \left( h \right)=\frac{\Gamma \left( h \right)\Gamma \left( {{\text{ }\!\!\delta\!\!\text{ }}_{1}}+{{\text{ }\!\!\delta\!\!\text{ }}_{2}}+h \right)}{\Gamma \left( {{\text{ }\!\!\delta\!\!\text{ }}_{1}}+h \right)\Gamma \left( {{\text{ }\!\!\delta\!\!\text{ }}_{2}}+h \right)},\,\,\Delta \left( h \right)=\frac{\Gamma \left( {{\text{ }\!\!\delta\!\!\text{ }}_{1}}+h \right)\Gamma \left( {{\text{ }\!\!\delta\!\!\text{ }}_{2}}+h \right)}{\Gamma \left( h \right)\Gamma \left( {{\text{ }\!\!\delta\!\!\text{ }}_{1}}+{{\text{ }\!\!\delta\!\!\text{ }}_{2}}+h \right)},
\end{equation}
where  ${{\text{ }\!\!\delta\!\!\text{ }}_{1}}={{x}_{1}}\displaystyle\frac{\partial }{\partial {{x}_{1}}}$, ${{\text{ }\!\!\delta\!\!\text{ }}_{2}}={{x}_{2}}\displaystyle\frac{\partial }{\partial {{x}_{2}}}$.

Recently, Hasanov and Srivastava \cite{{HS06},{HS07}} generalized the operators $\nabla \left( h \right)$ and $\Delta \left( h \right)$ defined by (\ref{eq31}) in the forms
\begin{equation} \label{eq32}
{{\tilde{\nabla }}_{{{x}_{1}};{{x}_{2}},...,{{x}_{n}}}}\left( h \right)=\frac{\text{ }\!\!\Gamma\!\!\text{ }\left( h \right)\text{ }\!\!\Gamma\!\!\text{ }\left( {{\text{ }\!\!\delta\!\!\text{ }}_{1}}+...+{{\text{ }\!\!\delta\!\!\text{ }}_{n}}+h \right)}{\text{ }\!\!\Gamma\!\!\text{ }\left( {{\text{ }\!\!\delta\!\!\text{ }}_{1}}+h \right)\text{ }\!\!\Gamma\!\!\text{ }\left( {{\text{ }\!\!\delta\!\!\text{ }}_{2}}+...+{{\text{ }\!\!\delta\!\!\text{ }}_{n}}+h \right)},
\end{equation}
\begin{equation} \label{eq33}
{{\tilde{\Delta }}_{{{x}_{1}};{{x}_{2}},...,{{x}_{m}}}}\left( h \right)=\frac{\text{ }\!\!\Gamma\!\!\text{ }\left( {{\text{ }\!\!\delta\!\!\text{ }}_{1}}+h \right)\text{ }\!\!\Gamma\!\!\text{ }\left( {{\text{ }\!\!\delta\!\!\text{ }}_{2}}+...+{{\text{ }\!\!\delta\!\!\text{ }}_{n}}+h \right)}{\text{ }\!\!\Gamma\!\!\text{ }\left( h \right)\text{ }\!\!\Gamma\!\!\text{ }\left( {{\text{ }\!\!\delta\!\!\text{ }}_{1}}+...+{{\text{ }\!\!\delta\!\!\text{ }}_{n}}+h \right)},
\end{equation}
where ${{\text{ }\!\!\delta\!\!\text{ }}_{i}}={{x}_{i}}\displaystyle\frac{\partial }{\partial {{x}_{i}}}\,\,\,\left( i=1,...,n \right)$ and they obtained very interesting results.  For example,  a hypergeometric function in $n$ variables $F_{A}^{\left( n \right)}$ defined by formula (\ref{eq12}) has the following decomposition formula [15]
\begin{equation*}
F_{A}^{(n)}\left( a,{{b}_{1}},...,{{b}_{n}};{{c}_{1}},...,{{c}_{n}};{{x}_{1}},...,{{x}_{n}} \right)
\end{equation*}
\begin{equation*}
=\sum\limits_{{{i}_{2}},...,{{i}_{n}}=0}^{\infty }{{}}\frac{{{\left( a \right)}_{{{i}_{2}}+...+{{i}_{n}}}}{{\left( {{b}_{1}} \right)}_{{{i}_{2}}+...+{{i}_{n}}}}{{\left( {{b}_{2}} \right)}_{{{i}_{2}}}}...{{\left( {{b}_{n}} \right)}_{{{i}_{n}}}}}{{{i}_{2}}!...{{i}_{n}}!{{\left( {{c}_{1}} \right)}_{{{i}_{2}}+...+{{i}_{n}}}}{{\left( {{c}_{2}} \right)}_{{{i}_{2}}}}...{{\left( {{c}_{n}} \right)}_{{{i}_{n}}}}}x_{1}^{{{i}_{2}}+...+{{i}_{n}}}x_{2}^{{{i}_{2}}}...x_{n}^{{{i}_{n}}}
\end{equation*}
\begin{equation*}
\cdot F\left( a+{{i}_{2}}+...+{{i}_{n}},{{b}_{1}}+{{i}_{2}}+...+{{i}_{n}};{{c}_{1}}+{{i}_{2}}+...+{{i}_{n}};{{x}_{1}} \right)
\end{equation*}
\begin{equation} \label{eq34}
\cdot F_{A}^{(n-1)}\left( a+{{i}_{2}}+...+{{i}_{n}},{{b}_{2}}+{{i}_{2}},...,{{b}_{n}}+{{i}_{n}};{{c}_{2}}+{{i}_{2}},...,{{c}_{n}}+{{i}_{n}};{{x}_{2}},...,{{x}_{n}} \right),n\in \mathbb{N}\backslash \left\{ 1 \right\}.
\end{equation}

However, due to the recurrence of the formula (\ref{eq34}), additional difficulties may arise in the applications of this expansion. Further study of the properties of operators (\ref{eq32}) and (\ref{eq33}) showed that formula (\ref{eq34}) can be reduced to a more convenient form \cite{Eturk}
\begin{equation*}
F_{A}^{(n)}\left( a,{{b}_{1}},...,{{b}_{n}};{{c}_{1}},...,{{c}_{n}};{{x }_{1}},...,{{x }_{n}} \right)=\sum\limits_{\underset{(2\le i\le j\le n)}{\mathop{{{s}_{i,j}}=0}}\,}^{\infty }{\frac{{{(a)}_{A(n,n)}}}{\underset{(2\le i\le j\le n)}{\mathop{{{s}_{2,2}}!{{s}_{2,3}}!\cdot \cdot \cdot {{s}_{i,j}}!\cdot \cdot \cdot {{s}_{n,n}}!}}\,}}
\end{equation*}
\begin{equation} \label{eq35}
\cdot \prod\limits_{l=1}^{n}\left[{\frac{{{({{b}_{l}})}_{B(l,n)}}}{{{({{c}_{l}})}_{B(l,n)}}}x _{l}^{B(l,n)}F\left( a+A(l,n),{{b}_{l}}+B(l,n);{{c}_{l}}+B(l,n);{{x }_{l}} \right)}\right],\, n \in\mathbb{N},
\end{equation}
where  $$A(l,n)=\sum\limits_{i=2}^{l+1}{\sum\limits_{j=i}^{n}{{{s}_{i,j}}}},\,\,B(l,n)=\sum\limits_{i=2}^{l}{{{s}_{i,l}}+}\sum\limits_{i=l+1}^{n}{{{s}_{l+1,i}}}.$$

Using  the well-known formula \cite[p.64, (22)]{E1953}
$$F\left(a,b;c;x \right)={{\left(1-x \right)}^{-b}}F\left( c-a,b;c;\frac{x}{x-1} \right),$$
we obtain
\begin{equation*}
   {F} _{A}^{(n)}\left( a,{{b}_{1}},...,{{b}_{n}};{{c}_{1}},...,{{c}_{n}};x_1,...,x_n \right)=\sum\limits_{\underset{(2\le p\le q\le n)}{\mathop{{{s}_{p,q}}=0}}\,}^{\infty }{\frac{{{(a)}_{A(n,n)}}}{{{s}_{2,2}}!\cdot \cdot \cdot {{s}_{n,n}}!}}\prod\limits_{l=1}^{n}\left[ \left(1-x_l\right)^{-b_l}\left(\frac{x_l}{1-x_l}\right)^{B(l,n)}\right]
\end{equation*}
\begin{equation} \label{eq351}
  \cdot \prod\limits_{l=1}^{n}\left[ \frac{{{({{b}_{l}})}_{B(l,n)}}}{{{({{c}_{l}})}_{B(l,n)}}}  \\
F\left(\begin{array}{*{20}c} {{c}_{l}}-a+B(l,n)-A(l,n),{{b}_{l}}+B(l,n);\\{{c}_{l}}+B(l,n);\\\end{array}\frac{x_l}{x_l-1} \right) \right]. \\
\end{equation}

It should be noted that the symbolic operators ${{\delta }_{1}}$ and ${{\delta }_{2}}$ in the one-dimensional case take the form $\delta :=xd/dx$  and such an operator is used in solving problems of the operational calculus \cite[p.26]{P}.

We now introduce here the other multivariable analogues of the Burchnall-Chaundy symbolic operators $\nabla \left( h \right)$   and  $\Delta \left( h \right)$  defined by (\ref{eq31}):
\begin{equation} \label{eq36}
\tilde{\nabla }_{x,y}^{\left( n,p \right)}\left( h \right):=\frac{\Gamma \left( h \right)\Gamma \left( h+{{\delta }_{1}}+...+{{\delta }_{n}}-{{\sigma }_{1}}-...-{{\sigma }_{p}} \right)}{\Gamma \left( h+{{\delta }_{1}}+...+{{\delta }_{n}} \right)\Gamma \left( h-{{\sigma }_{1}}-...-{{\sigma }_{p}} \right)}
\end{equation}
\begin{equation} \label{eq37}
=\sum\limits_{k=0}^{\infty }{\frac{{{\left( -{{\delta }_{1}}-...-{{\delta }_{n}} \right)}_{k}}{{\left( {{\sigma }_{1}}+...+{{\sigma }_{p}} \right)}_{k}}}{{{(h)}_{k}}k!}},
\end{equation}
\begin{equation} \label{eq38}
\tilde{\Delta }_{x,y}^{\left( n,p \right)}\left( h \right):=\frac{\Gamma \left( h+{{\delta }_{1}}+...+{{\delta }_{n}} \right)\Gamma \left( h-{{\sigma }_{1}}-...-{{\sigma }_{p}} \right)}{\Gamma \left( h \right)\Gamma \left( h+{{\delta }_{1}}+...+{{\delta }_{n}}-{{\sigma }_{1}}-...-{{\sigma }_{p}} \right)}
\end{equation}
\begin{equation*} \label{eq39}
=\sum\limits_{k=0}^{\infty }{\frac{{{\left( {{\delta }_{1}}+...+{{\delta }_{n}} \right)}_{k}}{{\left( -{{\sigma }_{1}}-...-{{\sigma }_{p}} \right)}_{k}}}{{{(1-h)}_{k}}k!}},
\end{equation*}
where
$$x:=(x_1,...,x_n),\,  y:=(y_1,...,y_p); \,{{\delta }_{i}}:={{x}_{i}}\frac{\partial }{\partial {{x}_{i}}},\,\,{{\sigma}_{j}}:={{y}_{j}}\frac{\partial }{\partial {{y}_{j}}},\,   i=\overline{1,n},\,\,j=\overline{1,p};\,\,n,p\in \mathbb{N}.$$

In addition, we consider operators which are equal to the Hasanov and Srivastava symbolic operators $\tilde{\nabla }\left( h \right)$ and $\tilde{\Delta }\left( h \right)$ defined by (\ref{eq32}) and (\ref{eq33}):
$$\tilde{\nabla }_{x,-}^{(n,0)}\left( h \right):={{\tilde{\nabla }}_{{{x}_{1}}:{{x}_{2}},...,{{x}_{n}}}}\left( h \right),      \,\,\tilde{\Delta }_{x,-}^{(n,0)}\left( h \right):={{\tilde{\Delta }}_{{{x}_{1}}:{{x}_{2}},...,{{x}_{n}}}}\left( h \right),\,\,n\in \mathbb{N};$$
$$\tilde{\nabla }_{-,y}^{(0,p)}\left( h \right):={{\tilde{\nabla }}_{-{{y}_{1}}:-{{y}_{2}},...,-{{y}_{p}}}}\left( h \right),\,\,\tilde{\Delta }_{-,y}^{(0,p)}\left( h \right):={{\tilde{\Delta }}_{-{{y}_{1}}:-{{y}_{2}},...,-{{y}_{p}}}}\left( h \right),\,\,p\in \mathbb{N}.$$
It is obvious that
$$\tilde{\nabla }_{x,-}^{(1,0)}\left( h \right)=\tilde{\Delta }_{x,-}^{(1,0)}\left( h \right)=\tilde{\nabla }_{-,y}^{(0,1)}\left( h \right)=\tilde{\Delta }_{-,y}^{(0,1)}\left( h \right)=1.$$

\textbf{Lemma 1}. Let  be $f:=f\left( x,y \right)$ function with variables $x:=(x_1,...,x_n)$  and $y:=(y_1,...,y_p)$ . Then following equalities hold true for any $n,\,p\in \mathbb{N}$:
\begin{equation} \label{eq310}
       {{\left( -\sum\limits_{i=1}^{n}{{{x}_{i}}\frac{\partial }{\partial {{x}_{i}}}} \right)}_{k}}f={{(-1)}^{k}}k!\sum\limits_{\begin{array}{*{20}c}
 {{i}_{1}}\ge 0,...{{i}_{n}}\ge 0 \\
 {{i}_{1}}+...+{{i}_{n}}=k
\end{array}}{\frac{x_{1}^{{{i}_{1}}}}{{{i}_{1}}!}...\frac{x_{n}^{{{i}_{n}}}}{{{i}_{n}}!}D_{x}^{{{i}_{1}}+...+{{i}_{n}}}f},\,\,\,k\in \mathbb{N};
\end{equation}
\begin{equation} \label{eq311}
       {{\left( \sum\limits_{j=1}^{p}{{{y}_{j}}\frac{\partial }{\partial {{y}_{j}}}} \right)}_{k}}f=k!\sum\limits_{l=1}^{k}{C_{k-1}^{l-1}}\sum\limits_{\begin{array}{*{20}c}
 {{j}_{1}}\ge 0,...{{j}_{p}}\ge 0 \\
 {{j}_{1}}+...+{{j}_{p}}=l
\end{array}}{\frac{y_{1}^{{{j}_{1}}}}{{{j}_{1}}!}...\frac{y_{p}^{{{j}_{p}}}}{{{j}_{p}}!}D_{y}^{{{j}_{1}}+...+{{j}_{p}}}f},\,\,k\in \mathbb{N}.
\end{equation}

The lemma 1 is proved by method of mathematical induction.

Now, we apply the symbolic operators $\tilde{\nabla }_{x,y}^{\left( n,p \right)}\left( h \right)$ and $\tilde{\Delta }_{x,y}^{\left( n,p \right)}\left( h \right)$ to the studying of properties of confluent hypergeometric function  ${\rm{H}} _{A}^{(n,p)}$  defined by (\ref{eq22}).

Using the formulas (\ref{eq36}) and (\ref{eq38}), we obtain
\begin{equation} \label{eq312}
{\rm{H}} _{A}^{(n,p)}\left( a,{{b}_{1}},...,{{b}_{n}};{{c}_{1}},...,{{c}_{n}};x,y \right)=\tilde{\nabla }_{x,y}^{\left( n,p \right)}\left( a \right)F_{A}^{(n)}\left( a,{{b}_{1}},...,{{b}_{n}};{{c}_{1}},...,{{c}_{n}};x \right){}_{0}{{F}_{p}}\left( 1-a;-y \right),                             \end{equation}
$$F_{A}^{(n)}\left( a,{{b}_{1}},...,{{b}_{n}};{{c}_{1}},...,{{c}_{n}};x \right){}_{0}{{F}_{p}}\left( 1-a;y \right)=\tilde{\Delta }_{x,y}^{\left( n,p \right)}\left( a \right){\rm{H}} _{A}^{(n,p)}\left( a,{{b}_{1}},...,{{b}_{n}};{{c}_{1}},...,{{c}_{n}};x,-y \right),$$                               where
$${}_{0}{{F}_{p}}\left( 1-a;y \right):=\sum\limits_{{{j}_{1}},...,{{j}_{p}}=0}^{\infty }{\frac{1}{{{\left( 1-a \right)}_{{{j}_{1}}+...+{{j}_{p}}}}}}\frac{y_{1}^{{{j}_{1}}}}{{{j}_{1}}!}...\frac{y_{p}^{{{j}_{p}}}}{{{j}_{p}}!}.$$

Now considering the equalities (\ref{eq37}), (\ref{eq310}) and (\ref{eq311}) from the formula (\ref{eq312}), we obtain the following decomposition formula for the confluent hypergeometric function ${\rm{H}} _{A}^{(n,p)}$:
\begin{equation*}
   {\rm{H}}_{A}^{(n,p)}\left( a,{{b}_{1}},...,{{b}_{n}};{{c}_{1}},...,{{c}_{n}};x,y \right)=F_{A}^{(n)}\left( a,{{b}_{1}},...,{{b}_{n}};{{c}_{1}},...,{{c}_{n}};x \right){}_{0}{{F}_{p}}\left( 1-a;-y \right) \\
\end{equation*}
\begin{equation*}
 +\sum\limits_{k=1}^{\infty }{\sum\limits_{l=1}^{k}{\sum\limits_{\begin{array}{*{20}c}
 {{i}_{1}}\ge 0,...,{{i}_{n}}\ge 0,{{j}_{1}}\ge 0,...,{{j}_{p}}\ge 0, \\
 {{i}_{1}}+...+{{i}_{n}}=k,{{j}_{1}}+...{{j}_{p}}=l
\end{array}}{C_{k-1}^{l-1}\frac{{{\left( -1 \right)}^{k+l}}k!{{\left( {{b}_{1}} \right)}_{{{i}_{1}}}}...{{\left( {{b}_{n}} \right)}_{{{i}_{n}}}}}{{{\left( 1-a \right)}_{l}}{{\left( {{c}_{1}} \right)}_{{{i}_{1}}}}...{{\left( {{c}_{n}} \right)}_{{{i}_{n}}}}}\frac{x_{1}^{{{i}_{1}}}}{{{i}_{1}}!}...\frac{x_{n}^{{{i}_{n}}}}{{{i}_{n}}!}\frac{y_{1}^{{{j}_{1}}}}{{{j}_{1}}!}...\frac{y_{p}^{{{j}_{p}}}}{{{j}_{p}}!}}}} \\
\end{equation*}
\begin{equation} \label{eq313}
 \cdot F_{A}^{(n)}\left( a+k,{{b}_{1}}+{{i}_{1}},...,{{b}_{n}}+{{i}_{n}};{{c}_{1}}+{{i}_{1}},...,{{c}_{n}}+{{i}_{n}};x \right){}_{0}{{F}_{p}}\left( 1-a+l;-y \right), \\
\end{equation}
where a Lauricella function $F^{(n)}_A$ has the expansion in the form (\ref{eq35}) or (\ref{eq351}).

Expansions (\ref{eq351}) and (\ref{eq313}) will be used for studying properties of the fundamental solutions.

\section{Confluent hypergeometric function  ${\rm{H}}_{A}^{(n,p)}$  in case $p=1.$}

Confluent hypergeometric function  ${\rm{H}}_{A}^{(n,p)}$  in case $p=1$ has the form
\begin{equation*}
{\rm{H}}_{A}^{(n,1)}\left( a,{{b}_{1}},...,{{b}_{n}};{{c}_{1}},...,{{c}_{n}};\xi ,{{\eta }_{1}} \right)=
{\rm{H}} _{A}^{(n,1)}\left[ \begin{array}{*{20}c}
   a,{{b}_{1}},...,{{b}_{n}};  \\
   {{c}_{1}},...,{{c}_{n}};  \\
\end{array}\text{ }\!\!\xi\!\!\text{ },{{\text{ }\!\!\eta\!\!\text{ }}_{1}} \right]
\end{equation*}
\begin{equation} \label{eq401}
=\sum\limits_{{{m}_{1}},...,{{m}_{n+1}}=0}^{\infty }{\frac{{{\left( a \right)}_{{{m}_{1}}+...+{{m}_{n}}-{{m}_{n+1}}}}{{\left( {{b}_{1}} \right)}_{{{m}_{1}}}}...{{\left( {{b}_{n}} \right)}_{{{m}_{n}}}}}{{{m}_{1}}!...{{m}_{n+1}}!{{\left( {{c}_{1}} \right)}_{{{m}_{1}}}}...{{\left( {{c}_{n}} \right)}_{{{m}_{n}}}}}\text{ }\!\!\xi\!\!\text{ }_{1}^{{{m}_{1}}}...\text{ }\!\!\xi\!\!\text{ }_{n}^{{{m}_{n}}}\eta _{1}^{{{m}_{n+1}}}},\,\left| {{\text{ }\!\!\xi\!\!\text{ }}_{1}} \right|+...+\,\left| {{\text{ }\!\!\xi\!\!\text{ }}_{n}} \right|<1,
\end{equation}

By virtue of the well-known formula \cite[p.52, (3)]{SM}
$$\sum\limits_{{{m}_{1}},...,{{m}_{p}}=0}^{\infty }{f\left( {{m}_{1}}+...+{{m}_{p}} \right)\frac{x_{1}^{{{m}_{1}}}}{{{m}_{1}}!}...\frac{x_{p}^{{{m}_{p}}}}{{{m}_{p}}!}}=\sum\limits_{m=0}^{\infty }{f\left( m \right)\frac{{{\left( {{x}_{1}}+...+{{x}_{p}} \right)}^{m}}}{m!}}$$
the following equality holds
$${\rm{H}} _{A}^{(n,p)}\left( ...;{{\xi }_{1}},...,{{\xi }_{n}}\text{,}{{\eta }_{1}},...,{{\eta }_{p}} \right)={\rm{H}} _{A}^{(n,1)}\left( ...;{{\xi }_{1}},...,{{\xi }_{n}}\text{,}{{\eta }_{1}}+...+{{\eta }_{p}} \right).$$

Therefore, it is sometimes enough to study a function ${\rm{H}}_{A}^{(n,1)}$ instead of function ${\rm{H}}_{A}^{(n,p)}.$

The function $\text{ }\!\!\omega\!\!\text{ }\left( \xi ,{{\eta }_{1}} \right)={\rm{H}} _{A}^{(n,1)}\left( a,{{b}_{1}},...,{{b}_{n}};{{c}_{1}},...,{{c}_{n}};{{\xi }_{1}},...,{{\xi }_{n}},{{\eta }_{1}} \right)$ defined by (\ref{eq401}) satisfies the following system of the hypergeometric equations
\begin{equation} \label{eq402}
\left\{ \begin{array}{*{20}c}
   {{\xi }_{i}}\left( 1-{{\xi }_{i}} \right){{\omega }_{{{\xi }_{i}}{{\xi }_{i}}}}-{{\xi }_{i}}\sum\limits_{j=1,j\ne i}^{n}{{{\xi }_{j}}{{\omega }_{{{\xi }_{i}}{{\xi }_{j}}}}}+{{\xi }_{i}}{{\eta }_{1}}{{\omega }_{{{\xi }_{i}}{{\eta }_{1}}}}+\left[ {{c}_{i}}-\left( a+{{b}_{i}}+1 \right){{\xi }_{i}} \right]{{\omega }_{{{\xi }_{i}}}}\\
   -{{b}_{i}}\sum\limits_{j=1,j\ne i}^{n}{{{\xi }_{j}}{{\omega }_{{{\xi }_{j}}}}}+{{b}_{i}}{{\eta }_{1}}{{\omega }_{{{\eta }_{1}}}}-a{{b}_{i}}\omega =0,\,\,\,i=\overline{1,n},  \\
   {{\eta }_{1}}{{\omega }_{{{\eta }_{1}}{{\eta }_{1}}}}-\sum\limits_{j=1}^{n}{{{\xi }_{j}}{{\omega }_{{{\xi }_{j}}}}}+\left( 1-a \right){{\omega }_{{{\eta }_{1}}}}+\omega =0,\,\,\,\,\,\,\,\,\,\,\,\,\,\,\,\,\,\,\,\,\,\,\,\,\,\,\,\,\,\,\,\,\,\,\,\,\,\,\,\,\,\,\,\,\,\,\,\,\,\,\,\,\,\,\,\,\,\,\,  \\
\end{array} \right.
\end{equation}
and all solutions of system (\ref{eq402}) are expressed by the formula:
\begin{equation}\label{eq403}
{{\omega }_{k}}\left( \xi ,\eta_1  \right)={{C}_{k}}\prod\limits_{i=1}^{k}{\left[ \xi _{i}^{1-{{c}_{i}}} \right]}\cdot {\rm{H}}_{A}^{\left( n,1 \right)}\left[\begin{array}{*{20}c} a,{{b}_{1}}+1-{{c}_{1}},...,{{b}_{k}}+1-{{c}_{k}},{{b}_{k+1}},...,{{b}_{n}};\\
2-{{c}_{1}},...,2-{{c}_{k}},{{c}_{k+1}},...,{{c}_{n}};\\
\end{array}\xi,\eta_1 \right],
\end{equation}
where  ${{C}_{k}}$ are arbitrary constants, $k=\overline{0,n}$.

In addition, the second statement of the lemma 1 is greatly simplified:
\begin{equation} \label{eq404}
{{\left( {{y}_{1}}\frac{\partial }{\partial {{y}_{1}}} \right)}_{k}}f=\sum\limits_{l=1}^{k}{\displaystyle\frac{k!(k-1)!y_{1}^{l}}{l!(l-1)!(k-l)!}\displaystyle\frac{{{\partial }^{l}}f(x,{{y}_{1}})}{\partial y_{1}^{l}}},\,\,k\in \mathbb{N}.
\end{equation}

Now, taking into account the symbolic equality (\ref{eq312}) at $p=1$ and using the formulas (\ref{eq310}) and (\ref{eq404}) we get

\begin{equation*}
   {\rm{H}}_{A}^{(n,1)}\left( a,{{b}_{1}},...,{{b}_{n}};{{c}_{1}},...,{{c}_{n}};x,{{y}_{1}} \right)=F_{A}^{(n)}\left( a,{{b}_{1}},...,{{b}_{n}};{{c}_{1}},...,{{c}_{n}};x \right) {}_0F_1(1-a;-y_1) \\
\end{equation*}
\begin{equation*}
  +\sum\limits_{k=1}^{\infty }{\sum\limits_{l=1}^{k}{\sum\limits_{\begin{array}{*{20}c}
 {{i}_{1}}\ge 0,...,{{i}_{n}}\ge 0, \\
 {{i}_{1}}+...+{{i}_{n}}=k
\end{array}}{C_{k-1}^{l-1}\frac{{{\left( -1 \right)}^{k+l}}k!{{\left( {{b}_{1}} \right)}_{{{i}_{1}}}}...{{\left( {{b}_{n}} \right)}_{{{i}_{n}}}}}{{{\left( 1-a \right)}_{l}}{{\left( {{c}_{1}} \right)}_{{{i}_{1}}}}...{{\left( {{c}_{n}} \right)}_{{{i}_{n}}}}}\frac{x_{1}^{{{i}_{1}}}}{{{i}_{1}}!}...\frac{x_{n}^{{{i}_{n}}}}{{{i}_{n}}!}\frac{y_{1}^{l}}{l!}}}} \\
\end{equation*}
\begin{equation} \label{eq405}
\cdot F_{A}^{(n)}\left( a+k,{{b}_{1}}+{{i}_{1}},...,{{b}_{n}}+{{i}_{n}};{{c}_{1}}+{{i}_{1}},...,{{c}_{n}}+{{i}_{n}};x \right){}_0F_1(1-a+l;-y_1), \\
\end{equation}
where the Lauricella function $F^{(n)}_A$ has the expansion in the form (\ref{eq35}) or  (\ref{eq351}).

Expansion (\ref{eq405}) will be used for studying properties of the fundamental solutions of the equation (\ref{eq11}).

\section{Fundamental solutions}

We consider equation (\ref{eq11}) in the domain $R_{m}^{n+}.$  Let $x:=\left( {{x}_{1}},...,{{x}_{m}} \right)$ be any point and ${{x}_{0}}:=\left( {{x}_{01}},...,{{x}_{0m}} \right)$  be any fixed point of $R_{m}^{n+}.$ We  search for a solution of equation (\ref{eq11}) as follows:
\begin{equation*} \label{eq501}
u(x,{{x}_{0}})=P(r)\omega (\xi \text{,}{{\eta }_{1}}),
\end{equation*}
where
\begin{equation} \label{eq5011}
P(r)={{\left( {{r}^{2}} \right)}^{-\alpha }}, \,\alpha =\sum\limits_{i=1}^{n}{{{\alpha }_{i}}}-1+\frac{m}{2};\,
\end{equation}
\begin{equation} \label{eq5012}
{{r}^{2}}=\sum\limits_{i=1}^{m}{{{({{x}_{i}}-{{x}_{0i}})}^{2}}},\,\,\,\,r_{j}^{2}={{({{x}_{j}}+{{x}_{0j}})}^{2}}+\sum\limits_{i=1,i\ne j}^{m}{{{({{x}_{i}}-{{x}_{0i}})}^{2}}}, \,j=\overline{1,n};
\end{equation}
\begin{equation*}
\xi :=\left( {{\xi }_{1}},...,{{\xi }_{n}} \right),\,\,{{\xi }_{j}}=\frac{{{r}^{2}}-r_{j}^{2}}{{{r}^{2}}},\,j=\overline{1,n};\,\,\, {{\eta }_{1}}=-\frac{1}{4}{{\lambda }^{2}}{{r}^{2}}.
\end{equation*}
We calculate all necessary derivatives and substitute them into equation (\ref{eq11}):
\begin{equation*}
\sum\limits_{m=1}^{n}{{{A}_{m}}\frac{{{\partial }^{2}}\text{ }\!\!\omega\!\!\text{ }}{\partial \text{ }\!\!\xi\!\!\text{ }_{m}^{2}}+{{A}_{n+1}}\frac{{{\partial }^{2}}\text{ }\!\!\omega\!\!\text{ }}{\partial \eta _{1}^{2}}+\sum\limits_{m=1}^{n}{\sum\limits_{k=m+1}^{n}{{{B}_{mk}}\frac{{{\partial }^{2}}\text{ }\!\!\omega\!\!\text{ }}{\partial {{\text{ }\!\!\xi\!\!\text{ }}_{m}}\partial {{\text{ }\!\!\xi\!\!\text{ }}_{k}}}}}}+\sum\limits_{m=1}^{n}{{{B}_{m,n+1}}\frac{{{\partial }^{2}}\text{ }\!\!\omega\!\!\text{ }}{\partial {{\xi }_{m}}\partial {{\eta }_{1}}}}
\end{equation*}
\begin{equation} \label{eq502}
+\sum\limits_{m=1}^{n}{{{D}_{m}}\frac{\partial \text{ }\!\!\omega\!\!\text{ }}{\partial {{\text{ }\!\!\xi\!\!\text{ }}_{m}}}}+{{D}_{n+1}}\frac{\partial \text{ }\!\!\omega\!\!\text{ }}{\partial {{\eta }_{1}}}+E\text{ }\!\!\omega\!\!\text{ }=0,
\end{equation}
where
$${{A}_{k}}=-\frac{4P(r)}{{{r}^{2}}}\frac{{{x}_{0k}}}{{{x}_{k}}}{{\xi }_{k}}\left( 1-{{\xi }_{k}} \right),\,   {{B}_{k,n+1}}=\frac{4P(r)}{{{r}^{2}}}\frac{{{x}_{0k}}}{{{x}_{k}}}{{\xi }_{k}}{{\eta }_{1}}+\frac{{{\lambda }^{2}}}{2}P(r){{\xi }_{k}},k=\overline{1,n}; $$
$${{B}_{kl}}=\frac{4P(r)}{{{r}^{2}}}\left( \frac{{{x}_{0k}}}{{{x}_{k}}}+\frac{{{x}_{0l}}}{{{x}_{l}}} \right){{\xi }_{k}}{{\xi }_{l}},\,k,l=\overline{1,n}, k\ne l;{{A}_{n+1}}={{\lambda }^{2}}P(r){{\eta }_{1}},$$ $${{D}_{k}}=-\frac{4P(r)}{{{r}^{2}}}\left\{ \left( 2{{\alpha }_{k}}-\alpha {{\xi }_{k}} \right)\frac{{{x}_{0k}}}{{{x}_{k}}}-{{\xi }_{k}}\sum\limits_{m=1}^{n}{\frac{{{x}_{0m}}}{{{x}_{m}}}{{\alpha }_{m}}} \right\},\,k=\overline{1,n};
$$ $${{D}_{n+1}}=\frac{4P(r)}{{{r}^{2}}}{{\eta }_{1}}\sum\limits_{m=1}^{n}{\frac{{{x}_{0m}}}{{{x}_{m}}}{{\alpha }_{m}}}+{{\lambda }^{2}}P(r)\alpha ,\,E=-{{\lambda }^{2}}P+\frac{4\alpha P(r)}{{{r}^{2}}}\sum\limits_{m=1}^{n}{\frac{{{x}_{0m}}}{{{x}_{m}}}{{\alpha }_{m}}}.$$

Using the given representations of coefficients we simplify equation (\ref{eq502}) and obtain the following system of equations:
\begin{equation} \label{eq503}
\left\{ \begin{array}{*{20}c}
    {{\xi }_{i}}\left( 1-{{\xi }_{i}} \right)\displaystyle\frac{{{\partial }^{2}}\omega }{\partial \xi _{i}^{2}}-{{\xi }_{i}}\sum\limits_{j=1,j\ne i}^{n}{{{\xi }_{j}}\displaystyle\frac{{{\partial }^{2}}\omega }{\partial {{\xi }_{i}}\partial {{\xi }_{j}}}}+{{\xi }_{i}}{{\eta }_{1}}\displaystyle\frac{{{\partial }^{2}}\omega }{\partial {{\xi }_{i}}\partial {{\eta }_{1}}}
  +\left[ 2{{\alpha }_{i}}-\left( \alpha +{{\alpha }_{i}}+1 \right){{\xi }_{i}} \right]\displaystyle\frac{\partial \omega }{\partial {{\xi }_{i}}}\\-{{\alpha }_{i}}\sum\limits_{j=1,j\ne i}^{n}{{{\xi }_{j}}\displaystyle\frac{\partial \omega }{\partial {{\xi }_{j}}}}+{{\alpha }_{i}}{{\eta }_{1}}\displaystyle\frac{\partial \omega }{\partial {{\eta }_{1}}}-\alpha {{\alpha }_{i}}\omega =0,\,i=\overline{1,n}, \\
  {{\eta }_{1}}\displaystyle\frac{{{\partial }^{2}}\omega }{\partial \eta _{1}^{2}}-\sum\limits_{j=1}^{n}{{{\xi }_{j}}\displaystyle\frac{{{\partial }^{2}}\omega }{\partial {{\xi }_{j}}\partial {{\eta }_{1}}}}+\left( 1-\alpha  \right)\displaystyle\frac{\partial \omega }{\partial {{\eta }_{1}}}+\omega =0,\,\,\,\,\,\,\,\,\,\,\,\,\,\,\,\,\,\,\,\,\,\,\,\,\,\,\,\,\,\,\,\,\,\,\,\,\,\,\,\,\,\,\,\,\,\,\,\,\,\,\,\,\,\,\,\,\,\,\,\,\,\,\,\,\,  \\
\end{array} \right.
\end{equation}

Now, using the solutions (\ref{eq403}) of the system of equations (\ref{eq402}), it is easy to determine the solutions of the system (\ref{eq503}), and substituting these solutions in (\ref{eq501}), we obtain the fundamental solutions of equation (\ref{eq11}) in the form
\begin{equation} \label{eq504}
{{q}_{k}}\left( x,{{x}_{0}} \right)={{\gamma }_{k}}\prod\limits_{i=1}^{k}{{{\left( {{x}_{i}}{{x}_{0i}} \right)}^{1-2{{\alpha }_{i}}}}\cdot }{{r}^{-2{{{\tilde{\alpha }}}_{k}}}}{\rm{H}}_{A}^{\left( n,1 \right)}\left[ \begin{array}{*{20}c}
   {{{\tilde{\alpha }}}_{k}},1-{{\alpha }_{1}},...,1-{{\alpha }_{k}},{{\alpha }_{k+1}},...,{{\alpha }_{n}};  \\
   2-2{{\alpha }_{1}},...,2-2{{\alpha }_{k}},2{{\alpha }_{k+1}},...,2{{\alpha }_{n}};  \\
\end{array}\xi ,\eta_1  \right],
\end{equation}
where
\begin{equation} \label{eq505}
{{\gamma }_{k}}={{2}^{2{{{\tilde{\alpha }}}_{k}}-m}}\frac{\Gamma \left( {{{\tilde{\alpha }}}_{k}} \right)}{{{\pi }^{m/2}}}\prod\limits_{j=1}^{k}{\frac{\Gamma \left( 1-{{\alpha }_{j}} \right)}{\Gamma \left( 2-2{{\alpha }_{j}} \right)}}\prod\limits_{i=k+1}^{n}{\frac{\Gamma \left( {{\alpha }_{i}} \right)}{\Gamma \left( 2{{\alpha }_{i}} \right)}},
\end{equation}
\begin{equation} \label{eq506}
{{\tilde{\alpha }}_{k}}=\frac{m}{2}+k-1-\sum\limits_{i=1}^{k}{{\alpha }_{i}}+\sum\limits_{i=k+1}^{n}{{\alpha }_{i}},\,k=\overline{0,n}.
\end{equation}

\section{Singularity properties of fundamental solutions}

Let us show that the found solutions (\ref{eq504}) have a singularity. We choose a solution
\begin{equation} \label{eq601}
{{q}_{0}}\left( x,{{x}_{0}} \right)={{\gamma }_{0}}{{r}^{-2{{{\tilde{\alpha }}}_{0}}}}{\rm{H}}_{A}^{\left( n,1 \right)}\left( {{{\tilde{\alpha }}}_{0}},{{\alpha }_{1}},...,{{\alpha }_{n}};2{{\alpha }_{1}},...,2{{\alpha }_{n}};\xi ,\eta_1  \right).
\end{equation}
For this aim we use the expansion (\ref{eq405}) for the confluent hypergeometric function ${\rm{H}}_{A}^{\left( n,1 \right)}\left( {{{\tilde{\alpha }}}_{0}},{{\alpha }_{1}},...,{{\alpha }_{n}};2{{\alpha }_{1}},...,2{{\alpha }_{n}};\xi ,\eta_1  \right)$.  As a result, solution (\ref{eq601}) can be written as follows
\begin{equation} \label{eq602}
{{q}_{0}}\left( x,{{x}_{0}} \right)={{r}^{2-m}}\prod\limits_{i=1}^{n}{\left[ r_{i}^{-2\alpha_i} \right]}\left\{ f_1\left( {{r}^{2}},r_{1}^{2},...,r_{n}^{2} \right){}_0F_1\left(1-a;{{\lambda }^{2}}{{r}^{2}}/4\right)+r^2f_2\left( {{r}^{2}},r_{1}^{2},...,r_{n}^{2};\lambda^2 \right)\right\},
\end{equation}
where
$$
 f_1\left( {{r}^{2}},r_{1}^{2},...,r_{n}^{2} \right)={\gamma}_{0}
 \sum\limits_{\underset{(2\le p\le q\le n)}{\mathop{{{s}_{p,q}}=0}}\,}^{\infty }{\frac{{{({{{\tilde{\alpha }}}_{0}})}_{A(n,n)}}}{{{s}_{2,2}}!\cdot \cdot \cdot {{s}_{n,n}}!}}\prod\limits_{l=1}^{n}\left[ \left(\frac{r^2}{r^2_l}-1\right)^{B(l,n)}\right]
$$
\begin{equation} \label{eq603}
  \cdot \prod\limits_{l=1}^{n}\left[ \frac{{{({{\alpha}_{l}})}_{B(l,n)}}}{{{({{2\alpha}_{l}})}_{B(l,n)}}}  \\
F\left(\begin{array}{*{20}c} {{2\alpha}_{l}}-{{{\tilde{\alpha }}}_{0}}+B(l,n)-A(l,n),{{\alpha}_{l}}+B(l,n);\\{{2\alpha}_{l}}+B(l,n);\\\end{array}1-\frac{r^2}{r^2_l} \right) \right] \\
\end{equation}
and $f_2\left( {{r}^{2}},r_{1}^{2},...,r_{n}^{2},\lambda^2 \right)$ has a finite value at $r\to 0 .$

It is easy to see that when $r\rightarrow 0$ the function $f_1\left( {{r}^{2}},r_{1}^{2},...,r_{n}^{2} \right)$ becomes an expression that does not depend on $x$ and $x_0$. Indeed, taking into account the equality
$$\prod\limits_{l=1}^{n}\left[ \left(-1\right)^{B(l,n)}\right]=1,$$
we have
\begin{equation*}
{\mathop {\lim} \limits_{r \to 0}}f_1\left( {{r}^{2}},r_{1}^{2},...,r_{n}^{2} \right)={\gamma}_{0}
 \sum\limits_{\underset{(2\le p\le q\le n)}{\mathop{{{s}_{p,q}}=0}}\,}^{\infty }{\frac{{{({{{\tilde{\alpha }}}_{0}})}_{A(n,n)}}}{{{s}_{2,2}}!\cdot \cdot \cdot {{s}_{n,n}}!}}
\end{equation*}
\begin{equation*} \label{eq604}
\cdot \prod\limits_{l=1}^{n}\left[ \frac{{{({{\alpha}_{l}})}_{B(l,n)}}}{{{({{2\alpha}_{l}})}_{B(l,n)}}}  \\
F\left(\begin{matrix} {{2\alpha}_{l}}-{{{\tilde{\alpha }}}_{0}}+B(l,n)-A(l,n),{{\alpha}_{l}}+B(l,n);\\{{2\alpha}_{l}}+B(l,n);\\\end{matrix}1 \right) \right] \\.
\end{equation*}

Applying now the summation formula \cite[p.104, (46)]{E1953}
$$
F(a,b;c;1)=\frac{\Gamma(c)\Gamma(c-a-b)}{\Gamma(c-a)\Gamma(c-b)}, \,\, Re(c-a-b)>0
$$
to each hypergeometric function $F(a,b;c;1)$ in the sum (\ref{eq603}), we get
\begin{equation*}
{\mathop {\lim} \limits_{r \to 0}}f_1\left( {{r}^{2}},r_{1}^{2},...,r_{n}^{2} \right)={\gamma}_{0}
 \prod\limits_{l=1}^{n}\left[\frac{\Gamma(2\alpha_l)\Gamma(\tilde{\alpha}_0-\alpha_l)}{\Gamma(\alpha_l)\Gamma(\tilde{\alpha}_0)}\right]
 \end{equation*}
\begin{equation*} \label{eq605}
\cdot\sum\limits_{\underset{(2\le p\le q\le n)}{\mathop{{{s}_{p,q}}=0}}\,}^{\infty }{\frac{{{({{{\tilde{\alpha }}}_{0}})}_{A(n,n)}}}{{{s}_{2,2}}!\cdot \cdot \cdot {{s}_{n,n}}!}} \prod\limits_{l=1}^{n}\left[ \frac{{{({{\alpha}_{l}})}_{B(l,n)}(\tilde{\alpha}_0-\alpha_l)_{A(l,n)-B(l,n)}}}{{{({\tilde{\alpha}_{0}})}_{A(l,n)}}}  \\
 \right] \\.
\end{equation*}

It is easy to verify that
$A\left( l,n \right)-B\left( l,n \right)
=\sum\limits_{i=2}^{l}{ \sum\limits_{j=i,j\neq l}^{n}{{{s}_{i,j}}}}
\ge 0,\,\,\,1\le l\le n.$

Taking into account the identity \cite[p.94,(33)]{Eturk}
$$
\sum\limits_{\underset{(2\le p\le q\le n)}{\mathop{{{s}_{p,q}}=0}}\,}^{\infty }{\frac{{{({{{\tilde{\alpha }}}_{0}})}_{A(n,n)}}}{{{s}_{2,2}}!\cdot \cdot \cdot {{s}_{n,n}}!}} \prod\limits_{l=1}^{n}\left[ \frac{{{({{\alpha}_{l}})}_{B(l,n)}(\tilde{\alpha}_0-\alpha_l)_{A(l,n)-B(l,n)}}}{{{({\tilde{\alpha}_{0}})}_{A(l,n)}}}  \\
 \right] =\Gamma\left(\frac{m}{2}\right)\frac{\Gamma^{n-1}({\tilde{\alpha}_{0}})}{\prod\limits\Gamma(\tilde{\alpha}_0-\alpha_l)}
$$
and
$$
{{\gamma }_{0}}={{2}^{2{{{\tilde{\alpha }}}_{0}}-m}}\frac{\Gamma \left( {{{\tilde{\alpha }}}_{0}} \right)}{{{\pi }^{m/2}}}\prod\limits_{i=1}^{n}{\frac{\Gamma \left( {{\alpha }_{i}} \right)}{\Gamma \left( 2{{\alpha }_{i}} \right)}},
$$
we obtain
\begin{equation} \label{eq606}
{\mathop {\lim} \limits_{r \to 0}}f_1\left( {{r}^{2}},r_{1}^{2},...,r_{n}^{2} \right)={{2}^{2{{{\tilde{\alpha }}}_{0}}-m}}\frac{\Gamma \left( m/2\right)}{{{\pi }^{m/2}}}.
\end{equation}
Expressions (\ref{eq602}) and (\ref{eq606}) give us the possibility to conclude that the solution $q_0(x,x_0)$ reduces to infinity of the order $r^{2-m}$ at $r\to 0$. Similarly, it is possible to be convinced that solutions $q_i(x,x_0), \, i=\overline{1,n},\,$ also reduce to infinity of the order $r^{2-m}$ when $r\to 0$.

Thus, all the fundamental solutions of the equation (\ref{eq11}) are written in the form (\ref{eq504}) and they have the singularity of the order $r^{2-m}\,(m>2)$ when $r\to 0$. We note that in the two-dimensional case (i.e. $m=2$) those solutions have a logarithmic singularity at  $r=0$ \cite{H07}.

\section{On fundamental solutions of the some generalized singular Helmholtz equation with "several" parameters}

 Consider the generalized Helmholtz equation with parameters  $\,\lambda_1,\lambda_2,...,\lambda_q\,$  before a required function
 \begin{equation} \label{eq71}
\sum\limits_{i=1}^{m}{\frac{{{\partial }^{2}}u}{\partial x_{i}^{2}}}+\sum\limits_{j=1}^{n}{\frac{2{{\alpha }_{j}}}{{{x}_{j}}}\frac{\partial u}{\partial {{x}_{j}}}}-\sum\limits_{k=1}^{p}{\lambda_k^2}u=0
\end{equation}
in the domain $R_{m}^{n+}:=\left\{ \left( {{x}_{1}},{{x}_{2}},...,{{x}_{m}} \right):{{x}_{1}}>0,{{x}_{2}}>0,...,{{x}_{n}}>0 \right\},$ where $m\ge 2$ is a dimension of the Euclidean space; $n\ge 1$ is a number of the singular coefficients;  $m\ge n;$ ${{\alpha }_{j}}$  are real constants and $0<2{{\alpha }_{j}}<1$$(j=\overline{1,n});$ $\lambda_k $ are real or pure imaginary constants, $k=\overline{1,p},\, p\in \mathbb{N}$.

Let $x:=\left( {{x}_{1}},...,{{x}_{m}} \right)$ be any point and ${{x}_{0}}:=\left( {{x}_{01}},...,{{x}_{0m}} \right)$  be any fixed point of $R_{m}^{n+}.$ We  search for a solution of equation (\ref{eq71}) as follows:
\begin{equation} \label{eq72}
u(x,{{x}_{0}})=P(r)\omega (\xi, \,{{\eta }}),
\end{equation}
where
\begin{equation*}
\xi :=\left( {{\xi }_{1}},...,{{\xi }_{n}} \right),\,\,\eta :=\left( {{\eta}_{1}},...,{{\eta }_{p}} \right); \,\,{{\xi }_{j}}=\frac{{{r}^{2}}-r_{j}^{2}}{{{r}^{2}}},\,j=\overline{1,n};\,\,\, {{\eta }_{k}}=-\frac{1}{4}{\lambda_k^{2}}{{r}^{2}},\,\, k=\overline{1,p};
\end{equation*}
the expressions  $P(r),\,r^2$ and $r_j^2$ are defined by formulas (\ref{eq5011}) and (\ref{eq5012}).

We calculate all necessary derivatives, substitute them into equation (\ref{eq11}) and obtain the following system of hypergeometric equations
\begin{equation} \label{eq702}
\left\{ \begin{array}{*{20}c}
   {{\xi }_{i}}\left( 1-{{\xi }_{i}} \right){{\omega }_{{{\xi }_{i}}{{\xi }_{i}}}}-{{\xi }_{i}}\sum\limits_{j=1,j\ne i}^{n}{{{\xi}_{j}}{{\omega }_{{{\xi }_{i}}{{\xi }_{j}}}}}+{{\xi }_{i}}\sum\limits_{j=1}^{p}{{{\eta }_{j}}{{\omega }_{{{\xi }_{i}}{{\eta}_{j}}}}}+\left[ {{2\alpha}_{i}}-\left( \tilde\alpha_0+{{\alpha}_{i}}+1 \right){{\xi }_{i}} \right]{{\omega }_{{{\xi }_{i}}}}\\
   -{{\alpha}_{i}}\sum\limits_{j=1,j\ne i}^{n}{{{\xi }_{j}}{{\omega }_{{{\xi }_{j}}}}}+{{\alpha}_{i}}\sum\limits_{j=1}^{p}{{{\eta }_{j}}{{\omega}_{{{\eta }_{j}}}}}-\tilde\alpha_0{{\alpha}_{i}}\omega =0,\,\,\,i=\overline{1,n},  \\
   \sum\limits_{l=1}^{p}{{{\eta }_{l}}{{\omega }_{{{\eta }_{l}}{{\eta }_{j}}}}}-\sum\limits_{l=1}^{n}{{{\xi }_{l}}{{\omega }_{{{\xi}_{l}}{{\eta }_{j}}}}}+\left( 1-\tilde\alpha_0 \right){{\omega }_{{{\eta }_{j}}}}+\omega =0,\,\,j=\overline{1,p}. \\
\end{array} \right.
\end{equation}

Now using the finding idea of all the solutions of system (\ref{eq24}) we get the solutions of system (\ref{eq702}) which have the form (\ref{eq25}) and substituting the found solutions in the formula (\ref{eq72}) we obtain the fundamental solutions of equation (\ref{eq71})
\begin{equation} \label{eq73}
{{q}_{k}}\left( x,{{x}_{0}} \right)={{\gamma }_{k}}\prod\limits_{i=1}^{k}{{{\left( {{x}_{i}}{{x}_{0i}} \right)}^{1-2{{\alpha }_{i}}}}\cdot }{{r}^{-2{{{\tilde{\alpha }}}_{k}}}}{\rm{H}}_{A}^{\left( n,p \right)}\left[ \begin{array}{*{20}c}
   {{{\tilde{\alpha }}}_{k}},1-{{\alpha }_{1}},...,1-{{\alpha }_{k}},{{\alpha }_{k+1}},...,{{\alpha }_{n}};  \\
   2-2{{\alpha }_{1}},...,2-2{{\alpha }_{k}},2{{\alpha }_{k+1}},...,2{{\alpha }_{n}};  \\
\end{array}\xi ,\eta  \right],
\end{equation}
where
${{\gamma }_{k}}$ and ${{\tilde{\alpha }}_{k}}\,(k=\overline{0,n})$ are defined by the formulas (\ref{eq505}) and (\ref{eq506}).

Our study shows that the fundamental solutions (\ref{eq73}) of equation (\ref{eq71}) have the same singularity as the solutions of the equation (\ref{eq11}) when $r\to 0$.

 \end{document}